\documentclass[11pt]{article}
\usepackage{amsmath,amssymb,amsfonts,amsthm,graphicx}
\usepackage{hyperref}
\newtheorem{theorem}{Theorem}[section]
\newtheorem{proposition}[theorem]{Proposition}

\newtheorem{lemma}[theorem]{Lemma}
\newtheorem{prob}[theorem]{Problem}
\newtheorem{remark}[theorem]{Remark}
\newtheorem{definition}[theorem]{Definition}
\numberwithin{equation}{section}
\title{A GAUSS-KUZMIN THEOREM FOR CONTINUED FRACTIONS ASSOCIATED WITH NON-POSITIVE INTEGER POWERS OF AN INTEGER $m \geq 2$}
\author{
    Dan Lascu\footnote{e-mail: lascudan@gmail.com.}\nonumber \\
    \emph{\small Mircea cel Batran Naval Academy, 1 Fulgerului, 900218 Constanta,
    Romania} 
    }
\begin{document}
\maketitle
\thispagestyle{empty}
\begin{abstract}
We consider a family $\{\tau_m:m\geq 2\}$ of interval maps introduced by Hei-Chi Chan \cite{Chan-2006}
as generalizations of the Gauss transformation.
For the continued fraction expansion arising from $\tau_m$,
we solve its Gauss-Kuzmin-type problem by applying the method of Rockett and Sz\"usz \cite{RS-1992}.
\end{abstract}

\section{Introduction}
H.-C. Chan considered some continued fraction expansions related to random Fibonacci-type sequences \cite{Chan-2005, Chan-2006}. In \cite{Chan-2005}, he studied the continued fraction expansions of a real number in the closed interval $[0,1]$ whose digits are differences of consecutive non-positive integer powers of $2$ and solve the corresponding Gauss-Kuzmin-L\'evy theorem. In fact, Chan has studied the transformation related to this new continued fraction expansion and the asymptotic behaviour of its distribution function. Giving a solution to the Gauss-Kuzmin-L\'evy problem, he showed in Theorem 1.1 in \cite{Chan-2005} that the convergence rate involved is $\mathcal{O}(q^n)$ as $n \rightarrow \infty$ with $0 < q < 1$.

The purpose of this paper is to prove a Gauss-Kuzmin type problem for the continued fraction expansions of real numbers in $[0,1]$ whose digits are differences of consecutive non-positive integer powers of an integer $m \geq 2$. In this section, we show our motivation and main theorems.

\subsection{Gauss' Problem and its progress}

One of the first and still one of the most important results in the metrical theory of continued fractions is so-called Gauss-Kuzmin theorem. 
Any irrational $0<x<1$ can be written as the infinite regular continued fraction 
\begin{equation}
x = \displaystyle \frac{1}{a_1+\displaystyle \frac{1}{a_2+\displaystyle \frac{1}{a_3+ \ddots}}} :=[a_1, a_2, a_3, \ldots],
\end{equation}
where $a_n \in \mathbb{N}_+ : = \left\{1, 2, 3, \ldots\right\}$.

Roughly speaking, the metrical theory (or, as called by Khintchine, the measure theory) of continued fraction expansions is about properties of the sequence $(a_n)_{n \in \mathbb{N}_+}$. It started on 25th October 1800, with a note by Gauss in his mathematical diary (entry 113) \cite{Brez}. 
Define the {\it regular continued fraction transformation} $\tau$ on the closed interval $I:=[0,1]$ by
 
\begin{equation}
\tau (x) = \left\{\begin{array}{lll} 
\displaystyle \frac{1}{x}-\left\lfloor \displaystyle \frac{1}{x} \right\rfloor & \hbox{if} & x \neq 0, \\ 
\\
0 & \hbox{if} & x = 0
\end{array} \right. \label{tau}
\end{equation}
where $\left\lfloor \cdot \right\rfloor$ denotes the floor (or entire) function. 
In modern notation, Gauss wrote that ``for very simple argument" we have
\begin{equation}
\lim_{n \rightarrow \infty} \lambda \left(\tau^n \leq x\right) = \frac{\log(1+x)}{\log2} \quad (x \in I) \label{01.2}
\end{equation}
where $\lambda$ denotes Lebesgue measure on $I$ and $\tau^{n}$ is the $n$-th iterate of $\tau$.

Nobody knows how Gauss found (\ref{01.2}), and his achievement is even more remarkable if we realize that modern probability theory and ergodic theory started almost a century later. In general, finding the invariant measure is a difficult task.

Twelve years later, in a letter dated 30 January 1812, Gauss wrote to Laplace that he did not succeed to solve satisfactorily ``a curious problem" and that his efforts ``were unfruitful". In modern notation, this problem is to estimate the error 
\begin{equation}
e_n(x) := \lambda \left(\tau^{-n}[0, x]\right) - \frac{\log(1+x)}{\log2} \quad (n \geq 1, \ x \in I).
\end{equation}
This has been called \textit{Gauss' Problem}. It received a first solution more than a century later, when R.O. Kuzmin \cite{Kuzmin-1928} showed in 1928 that 
\begin{equation}
e_n(x) = \mathcal{O}(q^{\sqrt{n}})
\end{equation}
as $n \rightarrow \infty$, uniformly in $x$ with some (unspecified) $0 < q < 1$. 
This has been called the \textit{Gauss-Kuzmin theorem} or the \textit{Kuzmin theorem}.

One year later, using a different method, Paul L\'evy \cite{Levy-1929} improved Kuzmin's result by showing that 
\begin{equation}
\left|e_n(x)\right| \leq q^n,
\end{equation}
$n \in \mathbb{N}_+$, $0 \leq x < 1$, with $q = 3.5 - 2\sqrt{2} = 0.67157...$. The Gauss-Kuzmin-L\'evy theorem is the first basic result in the rich metrical theory of continued fractions. 

By such a development, generalizations of these problems for non-regular continued fractions are also called as the \textit{Gauss-Kuzmin problem}.

\subsection{Chan's continued fraction expansions}

In this paper, we consider a generalization of the Gauss transformation and prove an analogous result. 
Especially, we will solve its Gauss-Kuzmin problem in Theorem \ref{Th.G-K}.

This transformation was studied in detail by Chan in \cite{Chan-2006} and Lascu in \cite{L-2013}.

Fix an integer $m \geq 2$. In \cite{Chan-2006}, Chan shows that any $x \in \left[0, 1\right)$ can be written as the form
\begin{equation}
x = \frac{m^{-a_1}}{\displaystyle 1+\frac{(m-1)m^{-a_2}}{\displaystyle 1 + \frac{(m-1)m^{-a_3}}{\displaystyle 1 + \ddots} }}:=[a_1, a_2, a_3, \ldots]_m, \label{3.1}
\end{equation}
where $a_n$'s are non-negative integers. 
Such $a_n$'s are also called {\it incomplete quotients (or continued fraction digits)} 
of $x$ with respect to the expansion in (\ref{3.1}) in this paper.

This continued fraction is treated as the following 
dynamical systems.
\begin{definition} \label{def.1}
Fix an integer $m \geq 2$.
\begin{enumerate}
\item[(i)]
The measure-theoretical dynamical system $(I,{\cal B}_{I},\tau_m)$ is defined as follows:
$I:=[0,1]$, 
$\mathcal{B}_I$ denotes the $\sigma$-algebra of all Borel subsets of $I$, 
and $\tau_{m}$ is the transformation 
\begin{equation}
\tau_m:I\to I;\quad 
\tau_{m}(x):=
\left\{
\begin{array}{ll}
{\displaystyle \frac{1}{m-1}\left(\frac{1}{m^{i}\,x}-1\right)} & 
{\displaystyle \mbox{if } x\in I_i,}\\
\\
0& \mbox{if }x=0
\end{array}
\right. \label{1.700}
\end{equation} 
where $I_i:=\{x\in I: m^{-(i+1)}< x\leq  m^{-i}\}$ for 
$i\in {\mathbb N}:=\{0,1,2,\ldots\}$.
\item[(ii)]
In addition to (i),
we write $(I,{\cal B}_{I},\gamma_m,\tau_m)$ as
$(I,{\cal B}_{I},\tau_m)$ with 
the following probability measure $\gamma_m$ on $(I,{\cal B}_{I})$:
\begin{equation}
\gamma_m (A) := k_m
\int_{A} \frac{dx}{\{(m-1)x+1\}\,\{(m-1)x+m\}}
\quad (A \in {\mathcal{B}}_I) \label{1.4''}
\end{equation}
where 
\begin{equation}
k_m := \frac{(m-1)^2}{\log \left\{m^2/(2m-1)\right\}}. \label{1.100'}
\end{equation}

\end{enumerate}
\end{definition}

Define the {\it quantized index map} $\eta_m:I\to {\mathbb N}$ by
\begin{equation}
\eta_m(x) := \left\{\begin{array}{ll} 
\lfloor -\log_m \,x \,\rfloor & \hbox{if }  x \neq 0, \\ 
\\
\infty & \hbox{if }  x = 0.
\end{array} \right. \label{1.4}
\end{equation}
By definition, $\eta_m(m^{-\alpha})=\lfloor \alpha \rfloor$.
By using $\tau_m$ and $\eta_m$, 
the sequence $(a_{n})_{n \in \mathbb N_+}$ in (\ref{3.1}) is obtained as follows:
\begin{equation}
a_n= \eta_m(\tau_m^{n-1}(x)) \quad (n \geq 1) \label{1.3}
\end{equation}
with $\tau_m^0 (x) = x$.
In this way, $\tau_m$ gives the algorithm of Chan's continued fraction expansion (\ref{3.1}).
\begin{proposition} \label{prop.3.3.1}
Let 
$(I,{\cal B}_{I},\gamma_m,\tau_m)$ 
be as in 
Definition \ref{def.1}(ii).
\begin{enumerate}
\item[(i)]
$(I,{\cal B}_{I},\gamma_m,\tau_m)$ is ergodic.
\item[(ii)]
The measure $\gamma_m$ is invariant under $\tau_m$,
that is, 
$\gamma_m (A) = \gamma_m (\tau^{-1}_m(A))$ 
for any $A \in {\mathcal{B}}_I$.
\end{enumerate}
\end{proposition}
\noindent \textbf{Proof.} 
See \cite{Chan-2006, L-2013}.
\hfill $\Box$\\

By Proposition \ref{prop.3.3.1}(ii),
$(I,{\cal B}_{I},\gamma_m,\tau_m)$ is a 
``dynamical system" in the sense of Definition 3.1.3 in \cite{BG-1997}.

\subsection{Known results and applications}

For Chan's continued fraction expansions,
we show known results and their applications in this subsection.

\subsubsection{Known results for $m=2$ case}

For $(I,{\cal B}_I,\gamma_m,\tau_m)$ in 
Definition \ref{def.1}(ii),
assume $m=2$, that is, we consider only 
$(I,{\cal B}_I,\gamma_2,\tau_2)$
in here. 

In \cite{Chan-2004}, 
Chan proved a Gauss-Kuzmin-L\'evy theorem for the transformation $\tau_2$. 
He showed that the convergence rate of the $n$-th distribution function 
of $\tau_2$ to its limit is $\mathcal{O}(q^n)$ 
as $n \rightarrow \infty$ with $q \leq 0.880555$ uniformly in $x$.

In \cite{Sebe-2005, Sebe-2007}, Sebe investigated the Perron-Frobenius operator
of $\tau_2$ by replacing a probability measure 
of the measurable space $(I,\mathcal {B}_I)$. 
Especially, she studied the Perron-Frobenius operator $J$ 
of $(I,{\cal B}_I,\gamma_2,\tau_2)$,
that is, $J$ is a unique operator on $L^1(I,\gamma_{2})$ satisfies
\begin{equation}
\int_A \{Jf\}(x)\,d \gamma_2(x)=\int_{\tau_{2}^{-1}(A)}f(x)\,d \gamma_2(x) 
\quad
(A \in {\cal B}_I,\,f\in L^1(I,\gamma_{2})).
\label{1.40}
\end{equation}
The asymptotic behavior of $J$ was shown 
by using well-known general results \cite{IG-2009,IK-2002}. 
By Wirsing-type approach \cite{Wirsing}, Sebe obtained 
a better estimate of the convergence rate involved \cite{Sebe-2005}. 
In fact, 
its upper and lower bounds of the convergence rate
were obtained as
$\mathcal{O}(w^n)$ and $\mathcal{O}(v^n)$, respectively 
when $n \rightarrow \infty$, with  $w < 0.209364308$ and $v > 0.206968896$
(\cite{Sebe-2005}, Theorem 4.3).
They provide a near-optimal solution to the Gauss-Kuzmin-L\'evy problem.

Furthermore, by restricting the Perron-Frobenius operator to the Banach space of functions $f: I \to {\mathbb{C}}$ of bounded variation, Iosifescu and Sebe \cite{IS-2006} proved that the exact optimal convergence rate of $\gamma_a (s^a_n \leq x)$ to $\gamma_2 ([0,x])$ is $\mathcal{O}(g^{2n})$ as $n \rightarrow \infty$ uniformly in $x$. Here $g$ is the inverse of the Golden ratio, i.e., we have 
\begin{equation}
g = \frac{2}{\sqrt{5} + 1}, \quad g^2 + g = 1, \quad g^2 = \displaystyle \frac{3 - \sqrt{5}}{2} = 0.38196 \ldots. 
\end{equation}
For $a \geq 0$, define the sequence 
$(s_{n,a})_{n\in{\mathbb{N}}}$ recursively 
by $s_{n,a} :=2^{-a_n}/(1 + s_{n-1,a})$, $n \in {\mathbb{N}_+}$, 
with $s_{0,a} := a$.
Then it is an $I\cup\{a\}$-valued Markov chain 
on $(I, \mathcal{B}_I, \gamma_a)$  where 
$\gamma_a$ is the probability measure on $(I,\mathcal{B}_I)$ defined as 
the following distribution function
\begin{equation}
\gamma_a ([0,x]) =: \frac{(a+2)x}{x+a+1} \quad ( x \in I, \ a \geq 0).
\end{equation}
For $J$ in (\ref{1.40}),
let $J^{'}$ denote its restriction on $L^{\infty}(I)\subset L^{1}(I,\gamma_2)$.
From Proposition 2.1.10 in \cite{IK-2002}, 
we see that
$J^{'}$ is the transition operator 
of the Markov chain $(s_{n,a})_{n \in {\mathbb{N}}_+}$ 
on $(I, \mathcal{B}_I, \gamma_a)$ for any $a \geq 0$.

\subsubsection{Known results for $m \geq 3$ case}

For $(I,{\cal B}_I,\gamma_m,\tau_m)$ in Definition \ref{def.1}(ii), recall the main results in \cite{L-2013, Sebe-2010}.

In \cite{L-2013}, Lascu proved a Gauss-Kuzmin theorem for the transformation $\tau_m$. In order to solve the problem,  he applied the theory of random systems with complete connections (RSCC) by Iosifescu \cite{IG-2009}.
We remind that a random system with complete connections is a quadruple
\begin{equation}
\left\{\left(I, {\mathcal{B}}_I\right), \left(\mathbb{N}_+, {\mathcal{P}}(\mathbb{N}_+)\right), u, P\right\}, \label{72}
\end{equation}
where $u: I \times \mathbb{N} \rightarrow I$, 
\begin{equation}
u(x,i)=u_{m,i}(x) = \frac{m^{-i}}{(m-1)x+1} \quad (x \in I) \label{3.63}
\end{equation}
and $P$ is the transition probability function from $\left(I,{\mathcal{B}}_I\right)$ to $\left(\mathbb{N}, {\mathcal{P}}(\mathbb{N}) \right)$ given by 
\begin{equation}
P(x,i) = P_{m,i}(x) = \frac{(m-1)m^{-(i+1)}(x+1)(x+m)}{(x+(m-1)m^{-i}+1)(x+(m-1)m^{-(i+1)}+1)}. \label{3.40}
\end{equation}
Also, the associated Markov operator of RSCCs (\ref{72}) is denote by $U_m$ and has the transition probability function 
\begin{equation}
Q_m(x, A) = \sum_{i \in W_m(x,A)}P_{m,i}(x) \quad (x \in I, \ A \in {\mathcal{B}}_I), 
\end{equation}
where $W_m(x,A) = \left\{i \in \mathbb{N}: u_{m,i}(x) \in A \right\}$.

Using the asymptotic and ergodic properties of operators associated with RSCC (\ref{72}), i.e., the ergodicity of RSCC, he obtained a convergence rate result for the Gauss-Kuzmin-type problem.

For more details about using RSCCs in solving Gauss-Kuzmin-L\'evy-type theorems, see \cite{IG-2009,IT-1969,Sebe-2001,Sebe-2002,Sebe-2005b, SL-2013}.

By Wirsing-type approach \cite{Wirsing} to the Perron-Frobenius operator of the associated transformation
under its invariant measure, Sebe \cite{Sebe-2010} studied the optimality of the convergence rate. 
Actually, Sebe obtained upper and lower bounds of the
convergence rate which provide a near-optimal solution to the Gauss-Kuzmin-L\'evy problem. 
In the case $m=3$, the upper and lower
bounds of the convergence rate were obtained as $O(w^n_3)$ and
$O(v^n_3)$ respectively when $n \rightarrow \infty$, with $v_3 > 0.262765464$ and $w_3 <
0.264687208$.

\subsubsection{Application to the asymptotic growth rate of a Fibonacci-type sequence} 
We explain an application of $\tau_m$ to a Fibonacci-type sequence here.
As is known, the {\it Fibonacci sequence} $(F_n)$ is 
recursively
defined as follows:
\begin{equation}
F_0 = F_1 = 1,\quad 
F_{n} = F_{n-1} + F_{n-2}\quad (n\geq 2).
\end{equation}
Equivalently, 
$(F_n)$ is also defined by {\it Binet's formula}
$F_n = (G^{n+1} - G^{-(n+1)})/\sqrt{5}$
for $n \geq 0$
where $G:=(1+\sqrt{5})/2$ is the Golden ratio.
By this formula, 
the asymptotic growth rate of $(F_n)$ is obtained as follows:
\begin{equation}
\lim_{n \rightarrow \infty} \frac{1}{n} \log F_n 
= \log \frac{1+\sqrt{5}}{2} = 0.4812\ldots
\end{equation}

A {\it random Fibonacci sequence} $(f_n)$ is defined 
as (with fixed $f_1$ and $f_2$) 
%
%
\begin{equation}
\label{eqn:fnalpha}
f_n = \alpha(n) f_{n-1} + \beta(n) f_{n-2},
\end{equation}
where $\alpha(n)$ and $\beta(n)$ are random coefficients.
For such $(f_{n})$,
the quest for its asymptotic growth rate is more difficult. 
We show two examples of random Fibonacci sequences as follows.

\begin{enumerate}
\item[(i)]
Define the random Fibonacci sequence $(f_n)$ as 
\begin{equation}
f_1 = f_2 = 1,\quad 
f_n = \pm f_{n-1} \pm f_{n-2}, \label{1.2400}
\end{equation}
where the signs in (\ref{1.2400})
are chosen independently 
and with equal probabilities.
Recently, Viswanath \cite{Viswanath} proved that 
its asymptotic growth rate is given as
\begin{equation}
\label{1.2500}
\lim_{n \rightarrow \infty} \frac{1}{n} \log |f_n |
= \log (1.13198824\ldots) = 0.12397559\ldots
\end{equation}
with probability $1$. 
\item[(ii)]
Fix an integer $m\geq 2$.
Define the random Fibonacci sequence $( f_n )$ as
\begin{equation}
f_{-1} = 0,\ f_0 = 1,\ a_0 = 0,\ 
f_{n}=m^{a_n} f_{n-1}+(m-1)m^{a_{n-1}}f_{n-2}, \label{1.3'}
\end{equation}
where $a_n$'s  are as in (\ref{3.1}). 
By using the ergodicity of 
$(I,{\cal B}_{I},\gamma_m,\tau_m)$ (Proposition \ref{prop.3.3.1}(i)),
Chan proved that
its asymptotic growth rate is given as follows \cite{Chan-2006}:
\begin{eqnarray}
\eta_m := \lim_{n \rightarrow \infty} \frac{1}{n} \log f_n\, 
&=& \,k_m \int^{1}_{0} \frac{-\log t}{\{(m-1)t + 1\}\,\{(m-1)t + m\}}\,dt  \nonumber \\
\, &\leq& \,k_m \frac{3m-1}{2m(2m-1)} \label{1.2600}
\end{eqnarray}
where $k_m$ is as in (\ref{1.100'}).
\end{enumerate}

\subsubsection{A Khintchin-type result and entropy} 

In probabilistic number theory, statistical limit theorems are established in problems involving ``almost independent'' random variables. The non-negative integers $a_n$, $n \in \mathbb{N}_+$, define random variables on the measure space $\left(I, \mathcal{B}_I, \textbf{P}\right)$, where $\textbf{P}$ is a probability measure on $I$. 

Continued fraction expansions of almost all irrational numbers are not periodic. Nevertheless, we readily reproduce another famous probabilistic result. It is the asymptotic value $\chi_m$ of the geometric mean of $m^{a_1}, m^{a_2}, \ldots, m^{a_n}$, i.e.,
\begin{equation}
\chi_m := \lim_{n \rightarrow \infty} \log \left( m^{a_1+a_2+\cdots+a_n} \right)^{1/n}, \label{1.30'}
\end{equation}
where $a_n$'s are given in (\ref{1.3}). This is a Khintchin-type result and we obtain 
\begin{eqnarray}
\chi_m &=& (\log m) \int_{0}^{1} a_1(x)\rho_m(x)dx \nonumber \\
&=& \frac{k_m \log m}{(m-1)^2} \sum_{n=0}^{\infty}\log \left( 1+\frac{(m-1)^3}{m^{n+2}+2(m-1)m+\frac{(m-1)^2}{m^n}} \right)^n \label{1.22''}
\end{eqnarray}
for almost all real numbers $x = [a_1(x), a_2(x), a_3(x), \ldots]_m \in (0, 1)$. As it can see, $\chi_m$ is a constant independent of the value of $x$. 

As it is well known, entropy is an important concept of information in physics, chemistry, and information theory \cite{PY}. The connection between entropy and the transmission of information was first studied by Shannon in \cite{Sh}. The entropy can be seen as a measure of randomness of the system, or the average information acquired under a single application of the underlying map. Entropy also plays an important role in ergodic theory. Thus in 1958 Kolmogorov \cite{Kol} imported Shannon's probabilistic notion of entropy into the theory of dynamical systems and showed how entropy can be used to tell whether two dynamical systems are non-conjugate. Like Birkhoff's ergodic theorem \cite{PY} the entropy is a fundamental result in ergodic theory. For a measure preserving transformation, its entropy is often defined by using partitions, but in 1964 Rohlin \cite{Roh} showed that the entropy of a $\mu$-measure preserving operator $T : [a, b] \rightarrow [a, b]$ is given by the beautiful formula
\begin{equation}
h(T) := \int^{b}_{a} \log\left|T' (x)\right| d\mu(x). \label{1.23}
\end{equation}
From Rohlin's formula it follows that the entropy of the operator $\tau_m$ in (\ref{1.700}) on the unit interval with respect to the measure $\gamma_m$ in (\ref{1.4''}) is given by
\begin{eqnarray}
h(\tau_m) &=& \int^{1}_{0} \log\left|\tau'_m (x)\right| \rho_m(x)dx = \int^{1}_{0} \log \left( \frac{m^{-a_1(x)}}{(m-1)x^2}\right)\rho_m(x)dx \nonumber \\
         &=& 2 \eta_m - \chi_m - \log(m-1), \label{1.23'}
\end{eqnarray}
where $a_1$, $\eta_m$ and $\chi_m$ are given in (\ref{1.3}), (\ref{1.2600}) and (\ref{1.30'}), respectively.

\subsection{Main theorem and its consequences}
\subsubsection{Main theorem}

We show our main theorems in this subsection. Fix an integer $m \geq 2$. Let $k_m$ be as in (\ref{1.100'}) and $(I, {\mathcal B}_{I})$ be as in Subsection 1.2. If $x$ has the expansion in (\ref{3.1}) and $\tau_m$ is as in $\mathrm(\ref{1.700})$, then it appears the question about the asymptotic distribution of $\tau^n_m$. If we know this, then the corresponding probability that $a_{n+1} = i$ is simply written as $\displaystyle prob \left(m^{-(i+1)} < \tau^n_m < m^{-i} \right)$. We will show that the event $\tau^n_m \leq x$ has the following asymptotic probability: 
\begin{equation}
\omega_m (x) = \frac{k_m}{(m-1)^2} \log \frac{m((m-1)x+1)}{(m-1)x+m} \quad ( x \in I). \label{1.15}
\end{equation}
This result allows us to say that the probability density function 
\begin{equation}
\rho_m(x) = \frac{k_m}{((m-1)x+1)((m-1)x+m)} \label{1.4.1}
\end{equation}
is invariant under $\tau_m$: if a random variable $X$ in the unit interval has density $\rho_m$, then so does $\tau_m$. 
The reason for this invariance is that for $0 \leq x < x+h \leq 1$, $\tau_m$ lies between $x$ and $x+h$ if and only if there exists $i \geq 1$, so that $X$ lies between $1 / (x+i+h)$ and $1 / (x+i)$. Thus
\begin{equation}
prob (x \leq \tau_m \leq x+h) = \sum_{i \in \mathbb{N}_+} prob \left(\frac{1}{x+i+h} \leq X \leq \frac{1}{x+i}\right).
\end{equation}
Taking the limits as $h \rightarrow \infty$ gives, for an arbitrary probability density function $f$ for $X$, the corresponding density $Gf$ for $\tau_m$ is given a.e. in $I$ by the equation 
\begin{equation}
G f(x) = \sum_{i \in \mathbb{N}} \frac{(m-1)m^{-i}}{((m-1)x+1)^2} f\left(\frac{m^{-i}}{(m-1)x+1}\right). 
\end{equation}
Clearly, the operator $G : L^1 \rightarrow L^1$ admits the density function $\rho_m$ as an eigenfunction corresponding to the eigenvalue $1$, i.e., $G \rho_m = \rho_m$. Here $L^1$ denotes the Banach space of all complex functions $f : I \rightarrow \mathbb{C}$ for which $\displaystyle \int_I \left|f\right| d\lambda < \infty$. 

The only eigenvalue of modulus $1$ of $G$ is $1$ and this eigenvalue is simple. 

From another perspective, the operator $\tau_m$ is an ergodic operator on the unit interval \cite{Chan-2006}, $\rho_m$ is the density of the invariant measure, and $G$ is called \textit{transfer operator} for $\tau_m$ \cite{L-2013}. The transfer operator $G$ has the same analytical expression as the Perron-Frobenius operator of $\tau_m$ under the Lebesgue measure \cite{L-2013}.

Our main result is the following theorem.

\begin{theorem} $\mathrm{(Gauss \mbox{-} Kuzmin \ Theorem)}$ \label{Th.G-K}
Let $\tau_m$ and $\omega_m$ be as in $\mathrm(\ref{1.700})$ and $\mathrm(\ref{1.15})$, respectively. When a non-atomic probability measure $\mu$ on $(I, \mathcal{B}_I)$ is given, define functions $F_{m,n}$ ($n \geq 0$) on $I$ by
\begin{equation}
\left\{
\begin{array}{ll}
F_{m,0} (x) := &\mu ([0,x]),\\
\\
F_{m,n} (x) :=& \mu (\tau_m^n \leq x)\quad(n\geq 1)\\
\end{array}
\right. \label{1.11'}
\end{equation}
for $x \in I$. 
Then there exists a constant $0<q_m<1$ such that $F_{m,n}$ is written as 
\begin{equation}
F_{m,n}(x) = \omega_m(x) + {\mathcal O}(q_m^n). \label{6}
\end{equation}
\end{theorem}

\begin{remark}
{\rm
\begin{enumerate}

\item[(i)]
From (\ref{6}), we see that 
\begin{equation}
\lim_{n\rightarrow\infty} F_{m,n}(x) = \gamma_m ([0,x]), \label{1.19}
\end{equation}
where $\gamma_m$ is the measure defined in (\ref{1.4''}). 
In fact, the Gauss-Kuzmin theorem estimates the error 
\begin{equation}
e_{m,n}(x) = e_{m,n}(x, \mu) = \mu (\tau_m^n \leq x) - \gamma_m ([0,x]) \quad (x \in I).
\end{equation}

\item[(ii)]
The solution of this problem implies that $(a_n)_{n \in \mathbb{N}}$ in (\ref{1.3}) is exponentially $\psi$-mixing under $\gamma_m$ (and under many other probability measures including $\lambda$) \cite{IG-2009, IK-2002}, that is, 
\begin{equation}
\left| \gamma_m (A_1 \cap A_2) - \gamma_m (A_1)\gamma_m (A_2) \right| \leq C q^n \gamma_m (A_1)\gamma_m (A_2) \quad (n \in \mathbb{N}_+) 
\end{equation}
for any $A_1 \in \sigma (a_1, \ldots, a_k)$ (the $\sigma$-algebra generated by the random variables $a_1, \ldots, a_k$), $A_2 \in \sigma (a_{n+k}, a_{n+k+1}, \ldots)$ and $k \in \mathbb{N}_+$, with suitable positive constants $q<1$ and $C$. 

In turn, $\psi$-mixing implies lots of limit theorems in both classical and functional versions. To form an idea of the results to be expected it is sufficient to look at the corresponding results for the regular continued fraction expansions \cite{IK-2002}.

\item[(iii)]
In (\ref{1.19}) we emphasized the probabilistic nature of Gauss's result. Khinchin \cite{Khin-1964} and Doeblin \cite{Doeblin} found new probabilistic results on the regular continued fraction transformation. These type of results were establish also for the transformation $\tau_m$ \cite{Chan-2006, L-2013}. These results establish, among other properties, that the map $\tau_m$ is ergodic (Proposition \ref{prop.3.3.1}(i)). Kuzmin's theorem may then be rephrased by saying that the convergence encountered in the mixing process (the ``approach to equilibrium'') is in fact exponential. If we define the linear operator $\Pi_1$ by
\begin{equation}
\Pi_1 f(x) = \rho_m(x) \int_{I} f d\lambda \quad (f \in L^1, \ x \in I)
\end{equation}
then there exists $ 0 < q_m < 1$ such that  
\begin{equation}
\left\|G^n - \Pi_1\right\|_{\mbox{v}} \leq \mathcal{O}(q^n_m) \quad (n \rightarrow \infty). 
\end{equation}
The norm $\left\|\cdot\right\|_{\mbox{v}}$ is defined by $ \left\|f\right\| _{\mbox{v}} = \left\|f\right\| + ``\mbox{total variation of}$ $f$'' \cite{IK-2002}. 

\end{enumerate}
}
\end{remark}

\begin{prob}
\label{prob:first}
{\rm
\begin{enumerate}
\item[(i)]
Solve the Gauss-Kuzmin-L\'evy problem of 
$\tau_m$ for $m \geq 3$.
For example, study the optimality of the convergence rate. Use the same strategy as in \cite{IS-2006}
\item[(ii)]
It is known that the Riemann zeta function is written
by using a kind of Mellin transformation 
of the Gauss transformation  $\tau$ in (\ref{tau})
as follows \cite{Vepstas}:
\begin{equation}
\label{eqn:zeta}
\zeta(s)=\frac{1}{s-1}-s\int_0^1 \tau(x)x^{s-1}\,dx
\quad (0<\Re (s)<1).
\end{equation}
This is derived by using the Euler-Maclaurin summation formula 
(\cite{TH}, p14) and the definition of $\tau$.
Then, by replacing $\tau$ with $\tau_m$ in (\ref{1.700}),
can we regard 
\begin{equation}
\label{eqn:zetam}
Z_m(s):=\frac{1}{s-1}-s\int_0^1 \tau_m(x)x^{s-1}\,dx
\end{equation}
as a new zeta function?
\end{enumerate}
}
\end{prob}

The paper is organised as follows. In Section 2, we prove Theorem \ref{Th.G-K}. In section 2.1, the necessary results used to prove the Gauss-Kuzmin theorem for the continued fractions presented in Section 1. The essential argument of the proof is the Gauss-Kuzmin-type equation. We will also give some results concerning the behavior of the derivative of $\{F_{m,n}\}$ in (\ref{1.11'}) which will allow us to complete the proof of the Theorem \ref{Th.G-K} in Section 2.2. 

\section{Proof of Theorem \ref{Th.G-K}} 
In this section, we will prove Theorem \ref{Th.G-K}. Fix an integer $m \geq 2$.

\subsection{Necessary lemmas} 

In this subsection, we show some lemmas. First, we show that $\{F_{m,n}\}$ in (\ref{1.11'}) satisfy a Gauss-Kuzmin-type equation.
\begin{lemma} \label{G-K.eq.}
For functions $\{F_{m,n}\}$ in $\mathrm{(\ref{1.11'})}$, the following Gauss-Kuzmin-type equation holds:
\begin{equation}
F_{m,n+1} (x) = \sum_{i \in \mathbb{N}}\left\{F_{m,n}\left(\alpha^i\right) - F_{m,n}\left(\frac{\alpha^i}{1+(m-1)x}\right)\right\} \label{2}
\end{equation}
for $x \in [0,1]$ and $n \in \mathbb{N}$ where $\alpha = 1/m$.
\end{lemma}
\noindent\textbf{Proof.}
Let $I:=[0,1]$, $I_{m,n}=\{x \in I: \tau^n_m(x) \leq x \}$ and $I_{m,n,i} = \left\{ x \in I_{m,n}: \ \alpha^i /(1+(m-1)x) < \tau^n_m(x) < \alpha^i \right\}$.

From (\ref{1.700}) and (\ref{1.3}), we see that
\begin{equation}
\tau_m^{n}(x) = \displaystyle\frac{m^{-a_{n+1}(x)}}{1+(m-1)\tau_m^{n+1}(x)} \quad (n \in \mathbb{N}_+). \label{2.3}
\end{equation}
From the definition of $I_{m,n,i}$ and (\ref{2.3}) it follows that for any $n \in \mathbb{N}$, $I_{m,n+1} = \bigcup_{i \in \mathbb{N}} I_{m,n,i}$. From this and using the $\sigma$-additivity of $\mu$, we have that
\begin{equation}
\mu(I_{m,n+1})=\mu\left(\displaystyle \bigcup_{i \in \mathbb{N}} I_{m,n,i}\right)=\sum_{i \in \mathbb{N}}\mu(I_{m,n,i}). \label{2.5}
\end{equation}
Then (\ref{2}) holds because $F_{m,n+1} (x)= \mu(I_{m,n+1})$ and 
\begin{equation}
\mu(I_{m,n,i}) = F_{m,n}\left(\alpha^i\right) - F_{m,n}\left(\frac{\alpha^i}{1+(m-1)x}\right). \label{2.6}
\end{equation}
\hfill $\Box$
\begin{remark}
\rm{
Assume that for some $p\in\mathbb{N}$, the derivative $F'_p$ exists everywhere in $I$ and is bounded. Then it is easy to see by induction that $F'_{m,p+n}$ exists and is bounded for all $n \in \mathbb{N}_+$. This allows us to differentiate (\ref{2}) term by term, obtaining 
\begin{equation}
F'_{m,n+1}(x) = \sum_{i \in \mathbb{N}} \frac{(m-1)\alpha^i}{(1+(m-1)x)^2}F'_{m,n}\left(\frac{\alpha^i}{1+(m-1)x}\right). \label{3}
\end{equation}
}
\end{remark}

\noindent We introduce functions $\{ f_{m,n} \}$ as follows:
\begin{equation}
f_{m,n}(x) := (1+(m-1)x)(m+(m-1)x)F'_{m,n}(x) \quad (x \in I, \ n \in \mathbb{N}). \label{2.8'}
\end{equation}
Then (\ref{3}) is 
\begin{equation}
f_{m,n+1}(x) = \sum_{i \in \mathbb{N}} P_m^i((m-1)x)f_{m,n}\left(\frac{\alpha^i}{1+(m-1)x}\right), \label{4}
\end{equation}
where $P^i_m(x)$ is given in (\ref{3.40}).

For $i \in \mathbb{N}$, define $\delta_i$ and $\beta^i_m(x)$ by
\begin{eqnarray}
\delta_i &:=& \alpha^{i} - \alpha^{2i} \label{3.63'} \\
\beta^i_m(x) &:=& \frac{(m-1)\delta_i}{\left((m-1)x + (m-1)\alpha^i+1\right)^2} \quad (x \in I). \label{3.63''}
\end{eqnarray}
\noindent Then we get
\begin{eqnarray}
P^i_m ((m-1)x) = (m-1) \left[\alpha^{i+1} + \frac{\delta_i}{(m-1)x+(m-1)\alpha^{i}+1} \right. \nonumber\\
\left.- \frac{\delta_{i+1}}{(m-1)x + (m-1)\alpha^{i+1}+1}\right]. \label{2.10}
\end{eqnarray}

\begin{lemma} \label{lema2.4}
For $\{f_{m,n}\}$ in $\mathrm{(\ref{2.8'})}$, define $M_n := \displaystyle \max_{x \in I}|f'_{m,n}(x)|$. Then
\begin{equation}
M_{n+1} \leq q_m \cdot M_n, \label{14}
\end{equation}
where
\begin{equation}
q_m := (m-1)^2(m^2+1) \sum_{i \in \mathbb{N}} \frac{1}{\left(m^{i+1}+m-1\right)^2}. \label{15}
\end{equation}
\end{lemma}
\noindent\textbf{Proof.}
We have 
\begin{equation}
\left(P_m^i((m-1)x)\right)' = (m-1) \left( \beta^{i+1}_m(x) - \beta^i_m(x) \right). \label{2.16}
\end{equation}
Now from (\ref{4}) and by calculus, we have 
\begin{equation}
f'_{m,n+1}(x) = (m-1)^2 \sum_{i \in \mathbb{N}} A_i f'_{m,n}(\theta_i) - (m-1) \sum_{i \in \mathbb{N}} B_i f'_{m,n}\left(u^i_m(x)\right), \label{10}
\end{equation}
where 
\begin{equation}
A_i := u^{i+1}_m(x) \beta^{i+1}_m(x), \quad B_i := P_m^i((m-1)x)\frac{\alpha^i}{((m-1)x+1)^2}
\end{equation}
and
\begin{equation}
u^{i+1}_m(x) < \theta_i < u^i_m(x).
\end{equation}

Now (\ref{10}) implies
\begin{equation}
M_{n+1}  \leq  M_n \cdot \max_{x \in I} \left| (m-1)^2 \sum_{i \in \mathbb{N}}A_i + (m-1) \sum_{i \in \mathbb{N}} B_i\right|. \label{11}
\end{equation}
We now must calculate the maximum value of the sums in this expression. 

First, we note that
\begin{equation}
A_i \leq \frac{\alpha^{2i+2}}{((m-1)\alpha^{i+1}+1)^2} \label{12}
\end{equation}
where we use $\delta_i = \alpha^{i} - \alpha^{2i}$ and $0 \leq x \leq 1$.

Next, observe that the function 
\begin{equation}
h_m(x) := \frac{\alpha^iP_m^i((m-1)x)}{((m-1)x+1)^2}
\end{equation}
is decreasing for $x \in I$ and $i \in \mathbb{N}$. Hence, $h_m(x) \leq h_m(0)$. This leads to
\begin{equation}
\frac{\alpha^iP_m^i((m-1)x)}{((m-1)x+1)^2} \leq \frac{(m-1)\alpha^{2i}}{((m-1)\alpha^{i+1}+1)^2}. \label{13}
\end{equation}
The relations (\ref{11}), (\ref{12}) and (\ref{13}) imply (\ref{14}) and (\ref{15}).
$\hfill \Box$

\subsection{Proof of Theorem \ref{Th.G-K}}

Introduce a function $R_{m,n}(x)$ such that
\begin{equation}
F_{m,n}(x) = \omega_m(x) + R_{m,n}\left(\omega_m(x)\right). \label{7}
\end{equation}

Because $F_{m,n}(0)=0$ and $F_{m,n}(1)=1$, we have $R_{m,n}(0)=R_{m,n}(1)=0$. To prove Theorem \ref{Th.G-K}, we have to show the existence of a constant $0<q_m<1$ such that 
\begin{equation}
R_{m,n}(x) = {\mathcal O}(q_m^n). \label{8}
\end{equation}

If we can show that $f_{m,n}(x)=k_m+{\mathcal O}(q_m^n)$, then its integration will show the equation (\ref{6}).

To demonstrate that $f_{m,n}(x)$ has this desired form, it suffices to prove the following lemma.
\begin{lemma} \label{lema2.3}
For any $x \in I$ and $n \in \mathbb{N}$, there exists a constant $0<q_m<1$ such that 
\begin{equation}
f'_{m,n}(x) = {\mathcal O}(q_m^n). \label{2.7}
\end{equation}
\end{lemma}
\noindent \textbf{Proof}.
Let $q_m$ be as in (\ref{15}). Using Lemma \ref{lema2.4}, to show (\ref{2.7}) it is enough to prove that $q_m < 1$.
To this end, for $i \geq 2$, observe that
\begin{equation}
\frac{1}{\left(m^{i+1}+m-1\right)^2} \leq \frac{1}{m^2(m-1)^2(m^2+1)} \left(\frac{1}{m}\right)^i.
\end{equation}
Therefore
\begin{eqnarray}
q_m \leq (m-1)^2(m^2+1) \qquad \qquad \qquad \qquad \qquad \qquad \qquad \qquad \qquad \qquad \quad \quad \quad \ \nonumber \\
\times \left\{ \frac{1}{(2m-1)^2} + \frac{1}{m^2+m-1} + \frac{1}{m^2(m-1)^2(m^2+1)} \sum_{i \geq 2}\left(\frac{1}{m}\right)^i  \right\} \quad \ \nonumber \\
= (m-1)^2(m^2+1) \left\{ \frac{1}{(2m-1)^2} + \frac{1}{m^2+m-1} + \frac{1}{m^3(m-1)^3(m^2+1)} \right\} \nonumber
\end{eqnarray}
\begin{equation}
 \leq 1, \qquad \qquad \qquad \qquad \qquad \qquad \qquad \qquad \qquad \qquad \qquad \qquad \quad 
\end{equation}
for any $m \in \mathbb{N}$, $m \geq 2$. 
$\hfill \Box$
\\



\end{document}